\documentclass[12pt]{report}
\usepackage{amssymb,amsmath,makeidx,verbatim}
\newcommand{\N}{\mathbb{N}}

\newcommand{\Z}{\mathbb{Z}}

\newcommand{\Q}{I\!\!\!Q}

\newcommand{\be}{\begin{enumerate}}
\newcommand{\ee}{\end{enumerate}}
\newcommand{\bq}{\begin{eqnarray*}}
\newcommand{\eq}{\end{eqnarray*}}

\begin{document}
\newcommand{\disp}{\displaystyle}
\thispagestyle{empty}
\begin{center}
\textsc{Non-gaussian $r-$congruences}\\
\ \\
{Olufemi O. OYADARE}\\
\ \\
Department of Mathematics,\\
Obafemi Awolowo University,\\
Ile-Ife, $220005,$ NIGERIA.\\
\text{E-mail: \textit{femi\_oya@yahoo.com}}\\
\end{center}
\baselineskip 1pt
\begin{quote}
\textbf{Abstract.} This paper is an \textit{invitation} to the study and use of general theory of non-gaussian $r-$congruences in the theory of numbers. In this work we classify the two kinds of $r-$congruences that exist (namely the \textit{trivial} and \textit{non-trivial} types) and establish their foundational properties as well as their algebra. Among other results, we show that the consideration of the \textit{non-trivial} $r-$ congruences entails the elucidation of a distinguished \textit{cyclic} subgroup of the \textit{permutation group,} while the \textit{trivial} $r-$congruences (chief among which is the well-known \textit{gaussian congruence}) leads to the trivial subgroup. The order of this cyclic group is also computed.
\end{quote}
\ \\
\ \\
\ \\
\ \\
\baselineskip 4pt

{\bf \S 1. Introduction.} The importance of congruence in the theory of numbers is evident from the fact that almost all (if not all) papers and texts on the subject, both at the elementary and advanced levels, contain and use results on it. In this paper we show that the well-known congruence, also called the \textit{gaussian congruence,} is one out of an infinite list of congruences, here called the \textit{$r-$congruences.} Our aim in this paper is to lay a foundation for the general theory of $r-$congruences thereby putting the gaussian congruence in its proper position among the other congruences.\\
\ \\

We believe that, with some efforts, results on gaussian congruence will ultimately be established for $r-$congruences while new connections will emerge from entirely new results on $r-$congruences.\\
\ \\
\ \\
\ \\
\ \\

{\bf \S 2. Main results.} Given $a, \;b \in \Z,$ with $a>0,$ there are unique integers $q$ and $r$ such that $b=aq+r$ where $0\leq r<a$ (or $\frac{-a}{2}<r\leq \frac{a}{2}$). There are two distinct ways of re-writing this equality in terms of the well-known gaussian congruence, $\equiv$: namely as $b\equiv r$ (mod $a$) (or as $b\equiv r$ (mod $q$), when $q\neq0$) and as $b\equiv aq$ (mod $r$) (if $r\neq0$). In this way we could not really write down any (gaussian) congruence between the initially given integers, $a$ and $b,$ except in the trivial case when $q=1$ (in which $a-b$ (if $a\neq b$) is thought of as the only divisor of itself). \textit{This way $7$ is not congruent to $8$ modulo $4.$} In this paper we seek to study a congruence that would always be possible between any given two integers, $a$ and $b.$\\
\ \\

To this end we observe, for any $a,\;b \in \Z,$ that $a-b \in \Z$ and that there always exist $m,\;r \in \Z,$ $m\neq0,$ such that $$m \mid (a-b)\;\;\mbox{with remainder $r$}\;\;\;\;([1.]).$$ We express this universal truth by writing $a\equiv_{r}b\;(\mbox{mod}\;m).$ It is clear that $a\equiv_{0}b\;(\mbox{mod}\;m)$ coincides with the aforementioned gaussian congruence, commonly written as $a\equiv b\;(\mbox{mod}\;m).$ Our aim in this paper is to consider the general $r-$congruence $a\equiv_{r}b\;(\mbox{mod}\;m)$ for any $r \in \Z$ and study the extent to which it could inherit properties of the gaussian congruence. We start with a formal definition of the main concept.\\
\ \\
\ \\

{\bf Definition 2.1.} \textit{Let $r \in \Z$ be fixed. An integer $a$ is said to be $r-$congruent to another integer $b$ modulo $m \in \Z \setminus \{0\}$ if $a-b=mq+r,$ with $q \in \Z.$}\\
\ \\

According to the division algorithm (quoted above) it follows that $0\leq r<m$ or $\frac{-m}{2}<r\leq \frac{m}{2}.$ We shall express the relationship in this definition by writing $$a\equiv_{r}b\;(\mbox{mod}\;m).$$ If $a\equiv_{0}b\;(\mbox{mod}\;m)$ we simply write the gaussian notation $a\equiv b\;(\mbox{mod}\;m).$ We can classify $r-$congruences into two types as described in the following.\\
\ \\
\ \\

{\bf Definition 2.2.} \textit{An $r-$congruence $a\equiv_{r}b\;(\mbox{mod}\;m)$ is said to be of a trivial type whenever $r \in m\Z$ and a non-trivial type whenever $r \notin m\Z.$}\\
\ \\

The following results gives a characterization of trivial $r-$congruences, which also sheds light on the importance and distinctiveness of the non-trivial type.\\
\ \\
\ \\

{\bf Lemma 2.3.} \textit{An $r-$congruence is of a trivial type if, and only if, it is gaussian.}\\
\ \\
\ \\

{\bf Proof.} Let $a\equiv_{r}b\;(\mbox{mod}\;m),$ then $a-b=mq+r,$ for $q \in \Z.$ If this $r-$congruence is trivial then $r=ms,$ $s \in \Z.$ Hence $a-b=mq+r=mq+ms=m(q+s).$ That is, $a-b=mt,$ where $t=q+s \in \Z.$ Therefore $a\equiv b\;(\mbox{mod}\;m),$ as required.\\
\ \\

Conversely, if $a\equiv b\;(\mbox{mod}\;m)$ then $a-b=md,$ where $d \in \Z.$ However every $d \in \Z$ may be written as $d=u+v,$ for some $u,\;v \in \Z.$ Hence $a-b=m(u+v)=mu+mv$ implying that $a\equiv_{mv}b\;(\mbox{mod}\;m).$ Setting $r=mv \in m\Z$ concludes the proof$.\;\;\Box$\\
\ \\

A direct consequence of the above Lemma is that \textit{every non-trivial $r-$congruence cannot be reduced to or deduced from a gaussian congruence.} That is, $a\equiv_{r}b\;(\mbox{mod}\;m),$ for $r\notin m\Z,$ has nothing to do with $a\equiv b\;(\mbox{mod}\;m).$ This Lemma justifies our efforts to study and understand $r-$congruences in general. The following lemma uplifts some central properties of a gaussian congruence to the status of an $r-$congruence.\\
\ \\
\ \\

{\bf Lemma 2.4.}\\
\ \\

\textit{$(i.)$ $a\equiv_{0}a\;(\mbox{mod}\;m)$ $\forall\;a \in \Z.$}\\
\ \\

\textit{$(ii.)$ $a\equiv_{r}b\;(\mbox{mod}\;m)$ if, and only if, $b\equiv_{(-r)}a\;(\mbox{mod}\;m).$}\\
\ \\

\textit{$(iii.)$ If $a\equiv_{r_{1}}b\;(\mbox{mod}\;m)$ and $b\equiv_{r_{2}}c\;(\mbox{mod}\;m),$ then $a\equiv_{(r_{1}+r_{2})}c\;(\mbox{mod}\;m).\;\;\Box$}\\
\ \\

In the light of Lemma $2.3,$ Lemma $2.4$ shows that there is an equivalence relation on $\Z$ from the trivial $r-$congruences. This is however a well-known fact in the theory of gaussian congruence. We also have the following general result for a fixed $r-$congruence.\\
\ \\
\ \\

{\bf Lemma 2.5.} \textit{Let $m,r \in \Z,\;m>0$ be fixed. Then the relation of $r-$congruence (mod $m$) is an equivalence relation on $\Z,$ in which a typical equivalence class, $\overline{a}_{m,r},$ for any $a \in \Z$ is given as $\overline{a}_{m,r}=\overline{a}_{m}+r,$ where $\overline{a}_{m}=$ equivalence class of $a$ in the gaussian congruence (mod $m$).}\\
\ \\
\ \\

{\bf Proof.} For any $a,\;b \in \Z$ define the relation $\sim$ on $\Z$ as $$a\sim b \iff a\equiv_{r}b\;(\mbox{mod}\;m),$$ for some $m \in \Z\setminus\{0\}.$ Then $a\sim a,$ $\forall\;a \in \Z,$ from $(2.4)(i.);$ $a\sim b$ implies $b\sim a$ $\forall\;a,\;b \in \Z,$ from $(2.4)(ii.);$ if $a\sim b$ and $b\sim c$ then $a\sim c$ $\forall\;a,\;b,\;c \in \Z,$ from $(2.4)(iii.).$ It is clear that $\overline{a}_{m,r}=\overline{a}_{m}+r.\;\;\Box$\\
\ \\

The formula $\overline{a}_{m,r}=\overline{a}_{m}+r,$ for the fixed $m,r \in \Z$ in Lemma $2.5$ shows that each of the (mod $m-$) equivalence classes of a fixed $r-$congruence is formed by shifting the corresponding (mod $m-$) equivalence classes of the gaussian congruence by the factor $r.$ In other words, each (mod $m-$) equivalence class of a fixed $r-$congruence is another (mod $m-$) equivalence class of the gaussian congruence. That is, if $a \in \Z,$ then $\overline{a}_{m,r}=\overline{a}_{m}+r=\overline{b}_{m},$ for some $b \in \Z.$ This close tie between an arbitrary $r-$congruence and the gaussian congruence may be seen from the fact that $$a\equiv_{r}b\;(\mbox{mod}\;m)\;\;\iff\;\;a\equiv (b+r)\;(\mbox{mod}\;m).$$\\
\ \\

This close tie does not however invalidate the importance and distinctiveness of the $r-$congruence and neither does it contradict Lemma $2.3$ above, since the $r-$congruence $a\equiv_{r}b\;(\mbox{mod}\;m),$ for $r\notin m\Z,$ has nothing to do with the gaussian congruence $a\equiv b\;(\mbox{mod}\;m).$ Indeed every $r-$congruence is in this way related to one another since $$a\equiv_{r_{1}}b\;(\mbox{mod}\;m)\;\;\iff\;\;a\equiv_{r_{2}} (b+r_{3})\;(\mbox{mod}\;m).$$ Here we have that $r_{1}=r_{2}+r_{3}.$ This general fact (among all the $r-$congruences) reduces to $$a\equiv_{r}b\;(\mbox{mod}\;m)\;\;\iff\;\;a\equiv (b+r)\;(\mbox{mod}\;m)$$ exactly in the special case where $r_{1}=r_{3}=:r.$ Hence there is nothing special about the gaussian congruence than that it was the first congruence to be discovered and that we are already so used to it (the same way we are used to \textit{denary numbers} at the detriment of other modes of numeration). The beauty of this relationship is in allowing us reduce the proofs of our results or calculations on $r-$congruence to a corresponding one in the gaussian congruence which will then be lifted back to the status of $r-$congruence at the end. The proofs of Lemmas $2.11$ and $2.12$ employ this technique in simplifying our arguments.\\
\ \\

Before discussing this further we consider an example of a non-trivial $r-$congruence in order to see the effect of the shift factor, $r,$ in the calculation of equivalence classes of $r-$congruences.\\
\ \\
\ \\

{\bf Example 2.6.} Here we give a computation for the equivalence classes of the $3-$congruence in modulo $5.$ Note in this case that $r=3 \notin 5\Z=m\Z,$ as required in $(2.2)$ and $(2.3)$ for non-trivial $r-$congruences. We have that $$\overline{0}_{5,3}=\overline{0}_{5}+3=(0+5\Z)+3=\{3+5n:n \in \Z\}=\{\cdots,-2,3,8,13,\cdots\}=\overline{3}_{5},$$ $$\overline{1}_{5,3}=\overline{1}_{5}+3=(1+5\Z)+3=\{4+5n:n \in \Z\}=\{\cdots,-1,4,9,14,\cdots\}=\overline{4}_{5},$$ $$\overline{2}_{5,3}=\overline{2}_{5}+3=(2+5\Z)+3=\{5+5n:n \in \Z\}=\{\cdots,-5,0,5,10,\cdots\}=\overline{0}_{5},$$ $$\overline{3}_{5,3}=\overline{3}_{5}+3=(3+5\Z)+3=\{6+5n:n \in \Z\}=\{\cdots,-4,1,6,11,\cdots\}=\overline{1}_{5}$$ and $$\overline{4}_{5,3}=\overline{4}_{5}+3=(4+5\Z)+3=\{7+5n:n \in \Z\}=\{\cdots,-3,2,7,12,\cdots\}=\overline{2}_{5}.$$ It is clear that if $\overline{a}_{m,r}=\overline{a}_{m}+r=\overline{b}_{m},$ for some $b \in \Z,$ then $b\equiv a+r\;(mod\;m).$\\
\ \\

We shall refer to each of $\overline{a}_{m,r},$ for a fixed $m,r \in \Z,\;m>0,$ as an \textit{$r-$residue class} of modulo $m.$\\
\ \\
\ \\

{\bf Lemma 2.7.} \textit{Let $\overline{a}_{m,r}$ denote an $r-$residue class of modulo $m.$ Then}\\
\ \\

\textit{$(i.)$ $\overline{a}_{m,r_{1}}=\overline{a}_{m,r_{2}}$ if, and only if, $r_{1}$ is a factor (or a multiple) of $r_{2}.$}\\
\ \\

\textit{$(ii.).$ $\overline{a}_{m,r}=\overline{b}_{m,r}$ if, and only if, $a\;(\mbox{mod}\;m)=b\;(\mbox{mod}\;m).\;\;\Box$}\\
\ \\

It then follows that $\overline{a}_{m,r+\beta}\neq\overline{a}_{m,r},$ for $\beta \in \Z$ in which $\beta\neq\alpha r,\;\forall\;\alpha \in \Z.$ \textit{The relationship between $\overline{a}_{m_{1},r}$ and $\overline{a}_{m_{2},r},$ for $m_{1},m_{2},r \in \Z\setminus\{0\},\;m_{1}\neq m_{2},$ is intricate and will be the topic of another paper.}\\
\ \\
\ \\

{\bf Theorem 2.8.} \textit{Let $r,\;m \in \Z,$ $m>0,$ be fixed. Then there are exactly $m$ distinct $r-$residue classes modulo $m.$}
\ \\
\ \\

{\bf Proof.} By Lemma $2.5,$ the $r-$residue classes modulo $m$ are $\overline{a}_{m,r}=\overline{a}_{m}+r=\overline{b}_{m},$ for some $b \in \Z.$ Since there are exactly $m$ distinct $0-$residue classes modulo $m$ and $r$ is fixed, then there are exactly $m$ distinct $r-$residue classes modulo $m.\;\;\Box$\\
\ \\

The above analysis of the $r-$residue classes of an $r-$congruence in relation to the gaussian congruence (as exemplified in Lemma $2.6$) reveals that the (mod $m-$) $r-$residue classes are derived by a \textit{re-arrangement} or \textit{permutation} of the (mod $m-$) gaussian residue classes. That this re-arrangement generates a distinguished \textit{cyclic} subgroup of the symmetric group, $S_{m},$ of order $m$ is established below.\\
\ \\
\ \\

{\bf Theorem 2.9.} \textit{Let $m,r \in \Z,$ $m>0,$ be fixed and let $\overline{a}_{m}$ denote the gaussian residue class of $a \in \Z$ modulo $m.$ Set\\
$$e_{m}=\left(\begin{array}{l}
\disp \overline{0}_{m}\;\;\;\;\overline{1}_{m}\;\;\;\;\overline{2}_{m}\;\;\;\;\cdots\;\;\;\;\overline{m-1}_{m}  \\
\disp \overline{0}_{m}\;\;\;\;\overline{1}_{m}\;\;\;\;\overline{2}_{m}\;\;\;\;\cdots\;\;\;\;\overline{m-1}_{m} \end{array}\right)\in S_{m}$$ and $$f^{(r)}_{m}=\left(\begin{array}{l}
\disp \;\;\;\overline{0}_{m}\;\;\;\;\;\;\;\;\;\;\overline{1}_{m}\;\;\;\;\;\;\;\;\;\;\overline{2}_{m}\;\;\;\;\;\;\;\cdots\;\;\;\;\;\;\overline{m-1}_{m} \\
\disp \overline{0}_{m}+r\;\;\;\;\overline{1}_{m}+r\;\;\;\;\overline{2}_{m}+r\;\;\;\;\cdots\;\;\;\;\overline{m-1}_{m}+r \end{array}\right)\;\in S_{m}.$$ Then the cyclic subgroup $C^{(r)}_{m}=\langle f^{(r)}_{m} \rangle$ of $S_{m}$ is trivial if $r\in m\Z$ and non-trivial of order $m$ if $r\notin m\Z.$}\\
\ \\

{\bf Proof.} It is clear that $f^{(r)}_{m}=e_{m}$ whenever $r \in m\Z,$ in which case $C^{(r)}_{m}$ is the trivial subgroup of $S_{m.}$ Now let $r \notin m\Z.$ We already know that $\overline{a}_{m,r}+s=\overline{a}_{m,r}$ iff $s\equiv0\;(mod\;m).$ That is, $\overline{a}_{m,r}+s=\overline{a}_{m,r}$ if and only if $s$ is a multiple of $m.$ Noting that $$f^{(r)}_{m}=\left(\begin{array}{l}
\disp \;\;\;\overline{0}_{m}\;\;\;\;\;\;\;\;\;\;\overline{1}_{m}\;\;\;\;\;\;\;\;\;\;\overline{2}_{m}\;\;\;\;\;\;\;\cdots\;\;\;\;\;\;\overline{m-1}_{m} \\
\disp \;\;\;\overline{0}_{m,r}\;\;\;\;\;\;\;\;\overline{1}_{m,r}\;\;\;\;\;\;\;\;\overline{2}_{m,r}\;\;\;\;\;\cdots\;\;\;\;\;\;\overline{m-1}_{m,r} \end{array}\right)\;\in S_{m},$$ it follows that $$\underbrace{f^{(r)}_{m}\circ f^{(r)}_{m}\circ \cdots \circ f^{(r)}_{m}}_{m-tuple}=e_{m},$$ which shows that the element $f^{(r)}_{m}$ of the symmetric group, $S_{m},$ has \textit{order} $m.$ Since the order of an element of a group is the same as the order of the cyclic (sub-)group it generates $([2.],\;p.\;59)$ we conclude that $C^{(r)}_{m}$ is non-trivial and is of order $m.\;\;\Box$\\
\ \\

The last Theorem gives credence to our choice of terms in Definition $2.2.$ The subgroup $C^{(r)}_{m}$ of $S_{m}$ (for $r \notin m\Z$), which shows the non-triviality of the concept of $r-$congruence, will be useful in understanding the group structure of non-trivial $r-$congruences and gives inter-connection among them. It should now be clear that number theory via $r-$congruences is a rich theory waiting to be explored. We consider next other properties of an $r-$congruence relating to its algebra.\\
\ \\
\ \\

{\bf Lemma 2.10.}\\
\ \\

\textit{$(i.)$ If $a\equiv_{r_{1}}b\;(\mbox{mod}\;m)$ and $c\equiv_{r_{2}}d\;(\mbox{mod}\;m)$ then $$a\pm c\equiv_{(r_{1}\pm r_{2})}b\pm d\;(\mbox{mod}\;m).$$ In particular, if $x \in \overline{a}_{m,r_{1}}$ and $\pm y \in \overline{b}_{m,r_{2}}$ then there is exactly one $r \in \Z$ for which $\pm x \in \overline{y}_{m,r}$ (or $\pm y \in \overline{x}_{m,r}$).}\\
\ \\

\textit{$(ii.)$ If $a\equiv_{r_{1}}b\;(\mbox{mod}\;m)$ and $c\equiv_{r_{2}}d\;(\mbox{mod}\;m)$ then $$ac\equiv_{(r_{1}d+r_{2}b+r_{1}r_{2})}bd\;(\mbox{mod}\;m).$$}\\
\ \\

\textit{$(iii.)$ If $a\equiv_{r}b\;(\mbox{mod}\;m)$ and $c \in \Z^{+}$ then $ca\equiv_{cr}cb\;(\mbox{mod}\;m).$}\\
\ \\

\textit{$(iv.)$ If $a\equiv_{r}b\;(\mbox{mod}\;m)$ and $d\mid m,$ where $d \in \Z^{+},$ then $a\equiv_{r}b\;(\mbox{mod}\;d).\;\;\Box$}\\
\ \\

\textit{$(v.)$ If $a\equiv_{r}b\;(\mbox{mod}\;m)$ and $k \in \N$ then $a^{k}\equiv_{f(r,b,k)}b^{k}\;(\mbox{mod}\;m),$ where the functions $(r,b,k)\mapsto f(r,b,k)$ is given as $f(r,b,k)=(r+b)^{k}-b^{k}.$}\\
\ \\
\ \\

{\bf Proof.} Items $(i.)-(iv.)$ are very clear while item $(v.)$ follows from mathematical inductions$.\;\;\Box$\\
\ \\

The well-known situation of Lemma $2.10$ for the gaussian congruence follows by setting $r=0.$ In this case $f(0,b,k)=0.$ In general, this Lemma reveals that the algebra of general $r-$congruence is dependent on $r,$ a situation totally absent for the gaussian congruence and that there are functions, like $(r,b,k)\mapsto f(r,b,k)$ in $2.10(v.)$ above, whose properties determine this algebra.\\
\ \\
\ \\

{\bf Lemma 2.11.} \textit{Let $m_{i} \in \Z^{+},$ $1\leq i\leq n,$ and denote the lcm of $m_{1},m_{2},\cdots,m_{n}$ by $[m_{1},m_{2},\cdots,m_{n}].$ Then $a\equiv_{r}b\;(\mbox{mod}\;m_{i})$ (for each $i\in\{1,2,\cdots,n\}$) if, and only if, $a\equiv_{r}b\;(\mbox{mod}\;[m_{1},m_{2},\cdots,m_{n}]).$}
\ \\
\ \\

{\bf Proof.} Let $a\equiv_{r}b\;(\mbox{mod}\;m_{i})$ (for each $i\in\{1,2,\cdots,n\}$). Since $a-b=m_{i}q_{i}+r,$ for some $q_{i},\;r \in \Z,$ then $a-b-r$ is a common multiple of all the $m_{i}'s.$ As $m_{i}\mid [m_{1},m_{2},\cdots,m_{n}]$ it follows that $[m_{1},m_{2},\cdots,m_{n}]\mid (a-b-r).$ That is, $a-b-r=[m_{1},m_{2},\cdots,m_{n}]q,$ for some $q \in \Z.$ Hence $a\equiv_{r}b\;(\mbox{mod}\;[m_{1},m_{2},\cdots,m_{n}]).$\\
\ \\

Conversely, let $a\equiv_{r}b\;(\mbox{mod}\;[m_{1},m_{2},\cdots,m_{n}]).$ Since $m_{i}\mid [m_{1},m_{2},\cdots,m_{n}],$ for each $i\in\{1,2,\cdots,n\},$ it follows that $a\equiv_{r}b\;(\mbox{mod}\;m_{i})$ (for each $i\in\{1,2,\cdots,n\}$)$.\;\;\Box$\\
\ \\

We have employed the fact that $$a\equiv_{r}b\;(\mbox{mod}\;m)\;\;\iff\;\;a\equiv (b+r)\;(\mbox{mod}\;m)$$ in the proof of the above result. One of the basic results needed to prove the \textit{Euler-Fermat Theorem} is the following.\\
\ \\
\ \\

{\bf Lemma 2.12.} \textit{Let $a,\;b,\;c \in \Z,$ $c\neq0,$ and let $ca\equiv_{r}cb\;(\mbox{mod}\;m),$ in which $c\mid r.$ Then $a\equiv_{r/c}b\;(\mbox{mod}\;m/(c,m)),$ where $(c,m)=$gcd of $c$ and $m.$}
\ \\
\ \\

{\bf Proof.} If $ca\equiv_{r}cb\;(\mbox{mod}\;m)$ then $ca\equiv (cb+r)\;(\mbox{mod}\;m)$ which implies that $a\equiv \frac{1}{c}(cb+r)\;(\mbox{mod}\;m/(c,m)).$ That is, $a\equiv (b+\frac{r}{c})\;(\mbox{mod}\;m/(c,m)).$ Hence $a\equiv_{r/c}b\;(\mbox{mod}\;m/(c,m)).\;\;\Box$\\
\ \\

A comparison of the classical ($r=0$) case of Lemma $2.12$ shows that the general situation above is only valid if $c \in \Z$ is chosen such that \textit{$c\mid r.$} That is, for $c \in \Z\setminus\{0\}$ in which $r \in c\Z.$ We know that every $c \in \Z\setminus\{0\}$ divides $0,$ which explains why the general requirement (that $c\mid r$) in the above Lemma is hidden in the classical case.\\
\ \\

Solvability of polynomial $r-$congruence may also be considered. We give a lemma in the case of linear $r-$congruences.\\
\ \\
\ \\

{\bf Lemma 2.13.} \textit{The linear $r-$congruence $ax\equiv_{r}b\;(\mbox{mod}\;m)$ is solvable for $x$ if, and only if, $(a,m)\mid(b+r).\;\;\Box$}\\
\ \\

It then follows that the notion of \textit{divisibility with zero-remainder} (commonly known as the gaussian congruence) is not the only concept of divisibility that exists or that should be studied in the theory of numbers. Indeed, as we have seen in this paper, divisibility with zero-remainder is only an example of a general concept of \textit{divisibility with $r-$remainder} (here called $r-$congruence), just as divisibility with $1-$remainder, divisibility with $2-$remainder, divisibility with $3-$remainder, ... are other examples worth understanding. Going back to the example mentioned at the start of this section, we can now conclude that:
\begin{quote}
\textit{It is unacceptable to say that $7$ is not congruent to $8$ modulo $4$ using only the restricted analysis of gaussian congruence, when in actual fact $$7\equiv_{r}8\;(mod\;4),\;\forall\;r \in 3+4\Z.$$}
\end{quote}
\ \\

{\bf Theorem 2.14.} \textit{Let $m \in \Z$ be fixed. The set $\mathfrak{R}_{m},$ of $r-$residue classes modulo $m,$ defined as $$\mathfrak{R}_{m}=\{\overline{a}_{m,r}:\;r \in \Z\}$$ and endowed with addition and multiplication as $$\overline{a}_{m,r_{1}}+\overline{a}_{m,r_{2}}:=\overline{a}_{m,r_{1}+r_{2}}\;\;\;\mbox{and}\;\;
\;\overline{a}_{m,r_{1}}\cdot\overline{a}_{m,r_{2}}:=\overline{a}_{m,r_{1}r_{2}}$$ is a ring, whose zero and identity elements are $\overline{a}_{m,0}$ and $\overline{a}_{m,1}.$}\\
\ \\
\ \\

{\bf Proof.} We define the map $\psi:\mathfrak{R}_{m}\rightarrow \Z$ as $\psi(\overline{a}_{m,r})=r$ and note that $\psi$ is a well-defined bijective ring homomorphism$.\;\Box$\\
\ \\

Analysis of the \textit{kernel} and \textit{co-kernel} of the ring homomorphism of Theorem $2.14$ may be a key to a complete understanding of the relationship between $\overline{a}_{m_{1},r}$ and $\overline{b}_{m_{2},r},$ for all $m_{1},m_{2},r \in \Z\setminus\{0\},\;m_{1}\neq m_{2}$ and $a,b \in \Z.$\\
\ \\

We believe that gaussian congruence is to $r-$congruence as Rolle's theorem (of Real analysis) is to the Mean-value theorem, in terms of its generality, and as the real number system is to the (infinitely many) $p-$adic number systems, in providing alternative view of the subject leading to a complete understanding of the notion of completion (of $\Q$). General $r-$congruence generalizes the gaussian congruence and provides a platform for the complete understanding of the notion of \textit{divisibility} in \textit{ideal} and \textit{module} theories.\\
\ \\

All the definitions and results of this paper and others on $r-$congruences may be stated for an arbitrary \textit{Euclidean domain} in place of $\Z.$\\
\ \\
\ \\
\ \\
\ \\
\ \\
\ \\
\ \\
\ \\
\ \\
\ \\
\ \\
\ \\
\ \\
\ \\

\indent {\bf References.}
\begin{description}
\item [{[1.]}] Cohn, H., \textit{Advanced Number Theory.} Dover Publications. $1980.$
\item [{[2.]}] Fraleigh, J. B., \textit{A First Course in Abstract Algebra.} Addison-Wesley, Reading. $2003.$
\end{description}
\end{document}